\documentclass[11pt]{amsart}
\usepackage{amsfonts,latexsym, color}
\usepackage{amsmath}
\usepackage{amssymb}
\usepackage[colorlinks=true, bookmarks=true, pdfstartview=FitH, pagebackref=true, linktocpage=true]{hyperref}

 \font\tenmsb=msbm10 at 12pt \font\sevenmsb=msbm7 at 8pt \font\fivemsb=msbm5 at
6pt
\newfam\msbfam
\textfont\msbfam=\tenmsb \scriptfont\msbfam=\sevenmsb \scriptscriptfont\msbfam=\fivemsb

\def\R{{\mathbb R}}

\def\N{{\mathbb N}}
\def\Sp{{\mathbb S}}

\begin{document}
\newcommand{\reset}{\setcounter{equation}{0}}

\newcommand{\beq}{\begin{equation}}
\newcommand{\noi}{\noindent}
\newcommand{\eeq}{\end{equation}}
\newcommand{\dis}{\displaystyle}
\newcommand{\mint}{-\!\!\!\!\!\!\int}

\def \theequation{\arabic{section}.\arabic{equation}}

\newtheorem{thm}{Theorem}[section]
\newtheorem{lem}[thm]{Lemma}
\newtheorem{cor}[thm]{Corollary}
\newtheorem{prop}[thm]{Proposition}
\theoremstyle{definition}
\newtheorem{defn}[thm]{Definition}
\newtheorem{rem}[thm]{Remark}

\def \bx{\hspace{2.5mm}\rule{2.5mm}{2.5mm}} \def \vs{\vspace*{0.2cm}} \def
\hs{\hspace*{0.6cm}}
\def \ds{\displaystyle}
\def \p{\partial}
\def \O{\Omega}
\def \H{{\mathbb H}}
\def \b{\beta}
\def \m{\mu}
\def \T{{\mathbb T}}
\def \ou{{\overline u}}
\def \ov{{\overline v}}
\def \D{\Delta}
\def \M{{\mathcal M}}
\def \l{\lambda}
\def \s{\sigma}
\def \e{\varepsilon}
\def \a{\alpha}
\def \o{\omega}
\def \b{\beta}
\def \SS{{\mathcal S}}
\def \Sp{{\mathbb S}}
\def \hw{{\widehat w}}
\def \HH{{\mathcal H}}
\def \H{{\mathbb H}}
\def \CC{{\mathcal C}}
\def \M{{\mathbb M}}
\def \E{{\mathbb E}}
\def\cqfd{%
\mbox{ }%
\nolinebreak%
\hfill%
\rule{2mm} {2mm}%
\medbreak%
\par%
}
\def \pr {\noindent {\bf Proof.} }
\def \rmk {\noindent {\it Remark} }
\def \esp {\hspace{4mm}}
\def \dsp {\hspace{2mm}}
\def \ssp {\hspace{1mm}}

\title{On Sharp Heisenberg Uncertainty Principle and the stability}
\author{Xia Huang}
\address{School of Mathematical Sciences, Key Laboratory of MEA (Ministry of Education) and Shanghai Key Laboratory of PMMP, East China Normal University, Shanghai, China}
\email{xhuang@cpde.ecnu.edu.cn}
\author{Dong Ye}
\address{School of Mathematical Sciences, Key Laboratory of MEA (Ministry of Education) and Shanghai Key Laboratory of PMMP, East China Normal
University, Shanghai, China}
\email{dye@math.ecnu.edu.cn}
\date{}
\begin{abstract}
In this work, we summarize the linearization method to study the Heisenberg Uncertainty Principles, and explain that the same approach can be used to handle the stability problem. As examples of application, combining with spherical harmonic decomposition and the Hardy inequalities, we revise two families of inequalities. We give firstly an affirmative answer in dimension four to Cazacu-Flynn-Lam's conjecture \cite{CFL} for the sharp Hydrogen Uncertainty Principle, and improve the recent estimates of Chen-Tang \cite{CT} in $\R^2$ and $\R^3$. On the other hand, we identify the best constants and extremal functions for two stability estimates associated to $\|\Delta u\|_2 \|r\nabla u\|_2 - \frac{N+2}{2}\|\nabla u\|^2_2$ in $\R^N$ ($N \geq 2$), studied recently by Duong-Nguyen \cite{DN} and Do-Lam-Lu-Zhang \cite{DLLZ}.
\end{abstract}

\maketitle

\section{Introduction}
\medskip
As a bridge, the Heisenberg Uncertainty Principle (HUP for shortness) translated a philosophical concept from quantum physics into a prolific research program in functional analysis. The core idea behind the famous HUP, that the non-commutativity of operators implies a lower bound for the measurement accuracy, played a foundational role in the study of functional inequalities. It yields not only various inequalities characterizing uncertainty, and presents also profound links with many important inequalities in analysis, as for example the Caffarelli-Kohn-Nirenberg inequalities.

Furthermore, it is well known that the optimal constant for the classical HUP is achieved by the Gaussian functions. This fact led mathematicians to systematically investigate the extremal functions and sharp constants for various uncertainty inequalities. The understanding of the sharpness of these functional inequalities provides interesting applications in the study of partial differential equations and harmonic analysis, see \cite{RS, F, FS, HJ, LS} and references therein.


\medskip
In general, the HUP can be seen as $L^2$ Caffarelli-Kohn-Nirenberg type inequalities, or $L^2$ interpolation inequalities as follows:
\begin{align}
\label{HUPg}
H(u)U(u) \geq \mu^2P(u)^2, \quad \forall\; u \in \H.
\end{align}
Here $\H$ is a functional Hilbert space, $H$, $U$, $P$ are continuous, positive definite quadratic forms over $\H$; and $\mu > 0$ is a constant. For simplicity, we consider only the real functions here, but there is no doubt to generalize to the complex value situations.

\medskip
Moreover, without loss of generality, we assume that $\H$ is invariant under the scaling of variables; while $H$, $U$, $P$ are homogeneous with respect to the scaling. To be more precise, let $\tau_\lambda$ stand for the homothety $\tau_\lambda(x) = \lambda x$ for $\lambda > 0$, and denote $u_\lambda(x) = u\circ \tau_\lambda(x) = u(\lambda x)$. We assume then $u_\lambda \in \H$ for any $\lambda > 0$ if $u \in \H$; and there exist $(\alpha_i) \in \R^3$ such that
\begin{align}
\label{scaling}
H(u_\lambda) = \lambda^{\alpha_1} H(u), \; U(u_\lambda) = \lambda^{\alpha_2} U(u), \; P(u_\lambda) = \lambda^{\alpha_3} P(u), \quad \forall\; u\in \H, \lambda > 0.
\end{align}
Obviously, to hope the inequality \eqref{HUPg}, we must have $\alpha_1 + \alpha_2 = 2\alpha_3$. Moreover, using the scaling, when $\gamma = \alpha_1 - \alpha_3 \ne 0$, it is easy to observe that \eqref{HUPg} is equivalent to the inequality:
\begin{align}
\label{HUP+}
H(u) + U(u) \geq 2\mu P(u), \quad \forall\; u \in \H.
\end{align}
The above fact is observed and used in \cite{CFL, CFLL, DN, H} for diverse situations, sometimes \eqref{HUP+} is called as {\it linearised} or {\it non scaling invariant} version of \eqref{HUPg}. Notice that conversely to the HUP, the inequality \eqref{HUP+} is no longer scaling invariant, but this enables an alternative way to handle \eqref{HUPg}, in particular associated to the spherical harmonic expansion.

Once we identify the best constant and extremal functions, we can ask for the stability of the Heisenberg uncertainty principle. Let $\Sigma_*$ denote the cone of extremal functions, i.e. functions such that the equalities holds in \eqref{HUPg} with the best constant, that is,
$$\mu_*^2 = \min_{u \in \H\backslash\{0\}} \frac{H(v)U(v)}{P(u)^2}\quad \mbox{and}\quad \Sigma_* = \Big\{ v \in \H, H(v)U(v) = \mu_*^2P(v)^2\Big\}\ne \{0_\H\}.$$
The question is to understand
\begin{align}
\label{staHUP}
\SS_* = \inf_{u \in \H, d_P(u, \Sigma_*)> 0} \frac{\sqrt{H(u)U(u)} - \mu_*P(u)}{d_P(u, \Sigma_*)^2}, \quad \mbox{where } \; d_P(u, \Sigma_*) = \inf_{v^* \in \Sigma_*} \sqrt{P(u - v^*)}.
\end{align}
Whenever $\SS_* > 0$, there holds
$$H(u)U(u) - \mu_*^2P(u)^2 \geq \SS_*^2 d_P(u, \Sigma_*)^4, \quad \forall\; u \in \H$$
which strengthen \eqref{HUPg}, by measuring quantitatively the deficit to the sharp estimate, in function of the distance to $\Sigma_*$, the set of extremal functions. The study of stability problem has also a {\it linearised} but equivalent form.
\begin{lem}
\label{newlem1}
Let $\H$, $H$, $U$, $P$, $\mu_*$ be as above, let $\SS_*$ be that in \eqref{staHUP}. Define
\begin{align}
\label{staHUP+}
\SS_{*, +} = \inf_{u \in \H, d_P(u, \Sigma_*)> 0} \frac{H(u) + U(u) - 2\mu_*P(u)}{d_P(u, \Sigma_*)^2}.
\end{align}
If $\alpha_1 \ne \alpha_3$ in \eqref{scaling}, then $\SS_{*, +} = 2\SS_*$.
\end{lem}

\medskip

The argument is the same as for the equivalence between \eqref{HUPg} and \eqref{HUP+}. Firstly, $H(u) + U(u)\geq 2\sqrt{H(u)U(u)}$ gives $\SS_{*, +} \geq 2\SS_*$. Inversely, remark that $d_P(u, \Sigma_*)^2$ has the same homogeneity as $P$ with respect to homotheties $\tau_\lambda$, that is, $d_P(u_\lambda, \Sigma_*)^2 = \lambda^{\alpha_3}d_P(u, \Sigma_*)^2$ for any $\lambda > 0$, due to the scaling invariance of $\Sigma_*$. Let $\gamma = \alpha_1-\alpha_3$, by definition of $\SS_{*, +}$, we get
\begin{align*}
\lambda^{\alpha_3} \SS_{*, +} d_P(u, \Sigma_*)^2 & = \SS_{*, +} d_P(u_\lambda, \Sigma_*)^2 \\
& \leq H(u_\lambda) + U(u_\lambda) - 2\mu_*P(u_\lambda)\\
& = \lambda^{\alpha_1} H(u) + \lambda^{\alpha_2}U(u) - 2\mu_*\lambda^{\alpha_3}P(u)\\
& = \lambda^{\alpha_3}\Big[\lambda^{\gamma} H(u) + \lambda^{-\gamma}U(u) - 2\mu_* P(u)\Big].
\end{align*}
Therefore $\lambda^{\gamma} H(u) + \lambda^{-\gamma}U(u) - 2\mu_* P(u) \geq \SS_{*, +} d_P(u, \Sigma_*)^2$ for any $\lambda > 0$. With minimization of the left hand side in $\lambda$, we arrive at
$$\sqrt{H(u)U(u)} - \mu_*P(u) \geq \frac{\SS_{*, +}}{2} d_P(u, \Sigma_*)^2.$$
It means $\SS_{*, +} \leq 2\SS_*$, hence $\SS_{*, +} = 2\SS_*$. 

Therefore, instead of studying \eqref{staHUP}, we can handle the {\it linearized} stability problem \eqref{staHUP+}.

\medskip
As applications, we revise two inequalities studied recently in \cite{CFL, CFLL, DLLZ, DN, DN1}. We will provide new sharpness results for the HUP or the corresponding stability estimate.

\medskip
A first example is the following Hydrogen Uncertainty Principle. More precisely, Cazacu-Flynn-Lam \cite{CFL} 
proved

\smallskip
\noindent
{\bf Theorem A.} {\sl Let $N\geq 5$ and $u\in W^{2,2}(\R^N)$. Then the following inequality holds
\begin{align}\label{CFL}
\int_{\R^N} |\Delta u|^2 dx \int_{\R^N}|\nabla u|^2 dx \geq \frac{(N+1)^2}{4}\left(\int_{\R^N}\frac{|\nabla u|^2}{|x|} dx\right)^2
\end{align}
where the constant $\frac{(N+1)^2}{4}$ is optimal, and is attained by $u(x)=\alpha(1+\beta r)e^{-\beta r},$ with any $\beta>0$, $\alpha\in \mathbb{R}$.}

\smallskip
In other words, let $\|\cdot\|_2$ denote the standard norm in $L^2(\R^N)$ and
\begin{align}\label{HN}
\HH_N:= \inf_{u\in W^{2,2}(\R^n), u \ne 0} \frac{\|\Delta u\|_2\|\nabla u\|_2}{\ds \int_{\R^N}\frac{|\nabla u|^2}{|x|}dx},
\end{align}
Cazacu-Flynn-Lam established that $\HH_N = \frac{N+1}{2}$ for $N \ge 5$. Moreover, they conjectured that \eqref{CFL} holds for $2\leq N\leq 4$. As the equality in \eqref{CFL} holds with $u(x) = \alpha(1+\beta r)e^{-\beta r}$, $\beta>0$ in any dimension $N \ge 2$, then \eqref{CFL} will be sharp whenever it is valid.

\medskip
Very recently, Chen-Tang give a negative answer to the above conjecture for dimensions two and three. More precisely, using explicit test functions, they obtain in \cite{CT} $$\frac{1}{2}\leq \HH_2\leq \frac{\sqrt{3}}{2} \quad \mbox{and}\quad \frac{3}{2}\leq \HH_3\leq \frac{2\sqrt{21}}{5} \simeq 1.833.$$
So the remained question is for $N =4$. Here we give a positive answer to \eqref{CFL} in $\R^4$, and present a unified study for $N\geq 4$.
\begin{thm}\label{1.1}
For $N\geq 4$, then $\HH_N=\frac{N+1}{2}$ and the extremal functions are
\begin{align*}
u(x)=\alpha(1+\beta r)e^{-\beta r},~\beta>0~ \text{and}~ \alpha\in \mathbb{R}.
\end{align*}
\end{thm}

For the inequality in \eqref{CFL}, we have $(\alpha_i) = (4-N, 2-N, 3-N)$. As explained abovely, to get the best constant $\HH_N^2$ for \eqref{CFL}, we can consider
\begin{align}
\label{HN+}
\HH_N^+ = 2\HH_N = \inf_{u\in W^{2,2}(\R^n), u \ne 0} \frac{\|\Delta u\|^2_2 + \|\nabla u\|^2_2}{\big\|\frac{1}{\sqrt{r}}\nabla u\big\|^2_2}.
\end{align}

Another ingredient of our study is based on the link between Hardy inequalities and HUP. A notable example is: Using the famous Hardy inequality
$$\|\nabla u\|_2 \geq \frac{N-2}{2}\Big\|\frac{u}{r}\Big\|_2, \quad \forall \; u \in W^{1,2}(\R^N).$$
we get immediately a Heisenberg type estimate via the Cauchy-Schwarz inequality
$$\|\nabla u\|_2\|ru\|_2 \geq \frac{N-2}{2}\|u\|_2^2, \quad \forall \; u \in W^{1,2}(\R^N).$$
Though we have not yet the best constant (which is $\frac{N}{2}$) in the above estimate, the Hardy inequality provides very useful information.

\smallskip
With the same approach, we improve also the estimates of Chen-Tang (see \cite[Theorem 1.1]{CT}) for $\HH_2$ and $\HH_3$. Another interesting fact is that without knowing the exact values of $\HH_2$ or $\HH_3$, we can claim the existence of extremal functions.
\begin{thm}
\label{HH2-3}
There hold
\begin{align}
\label{newHN2-3}
\frac{3+6\sqrt 2}{14}\leq \HH_2 \leq \frac{\sqrt 3}{2}, \quad 1.75\leq \HH_3\leq \frac{2\sqrt{21}}{5}.
\end{align}
Moreover, the best constant $\HH_2$ (resp. $\HH_3$) is attained in the closure of $C_c^\infty(\R^2)$ (resp. $C_c^\infty(\R^3)$) under the norm $\|v\| = \|\Delta v\|_2+\|\nabla v\|_2$.
\end{thm}

\medskip
As a second example of application, we consider the following HUP:
\begin{align}\label{HUP0}
\|\Delta u\|_2 \|r\nabla u\|_2 \geq \frac{N+2}{2}\|\nabla u\|^2_2.
\end{align}
Thanks to \cite{CFL} (see also \cite{DN1}), it is known that the constant $\frac{N+2}{2}$ in \eqref{HUP0} is optimal, and it is attained by the Gaussian functions
$$u \in \Sigma_0:= \Big\{\alpha e^{-\beta r^2}, \; \alpha\in \R, \; \beta>0 \Big\}.$$
The stability problem for \eqref{HUP0} was studied firstly by Duong-Nguyen \cite{DN}, they showed that for any $u \in \H_0$, 
\begin{align*}
\delta_1(u) := \|\Delta u\|_2 \|r\nabla u\|_2 - \frac{N+2}{2}\|\nabla u\|^2_2 \geq \frac{1}{768} \inf_{u^* \in S_0^u}\|\nabla(u-u^*)\|_2^2
\end{align*}
where
\begin{align*}
S_0^u:= \Big\{w\in \Sigma_0,\; \|\nabla w \|_2 =\|\nabla u\|_2\Big\}.
\end{align*}
Notice that $S_0^u$ is a sphere in $\Sigma_0$. Recently, applying the equivalent linearized approach \eqref{staHUP+}, Do-Lam-Lu-Zhang \cite{DLLZ} improved
drastically the estimate of Duong-Nguyen as follows.

\smallskip
\noindent
{\bf Theorem B.} {\sl Let $N \geq 2$, and $\H_0$ be the completion of $C_c^\infty(\R^N)$ under the norm $\|\D u\|_2 + \|r\nabla u\|_2$. Denote
$\CC_{N,1} = \sqrt{N^2 + 4N -4} -N$, then
\begin{align}
\label{staHUP0}
\delta_1(u) \geq \frac{\CC_{N,1}}{2}\inf_{u^* \in \Sigma_0}\|\nabla(u-u^*)\|_2^2, \quad \forall\; u \in \H_0
\end{align} and
\begin{align}
\label{staHUP0-DNbis}
\delta_1(u) \geq \frac{\CC_{N,1}}{4}\inf_{u^* \in S_0^u}\|\nabla(u-u^*)\|_2^2, \quad \forall\; u \in \H_0.
\end{align}
}

By refining the analysis in \cite{DLLZ}, we claim the optimality of the above inequalities.
\begin{thm}
\label{sharpsta}
Let $N\geq 2$, the above estimates \eqref{staHUP0} and \eqref{staHUP0-DNbis} are sharp. Moreover, the equalities occurs in \eqref{staHUP0} and \eqref{staHUP0-DNbis} for $u(x) = u_1(\lambda x)$, with $\lambda > 0$, $u(x) = r\Psi(r)\phi_1(\sigma)$ where
$$\Psi(r) :={_1F_1}\Big(\frac{N+1}{2} + \frac{\CC_{N,1}}{4};\frac{N+2}{2};-\frac{r^2}{2}\Big),$$
$\phi_1$ is any spherical eigenfunction on the unit sphere satisfying $-\D_{S^{N-1}}\phi_1 = (N-1)\phi_1$, and $\sigma = \frac{x}{|x|}$.
\end{thm}
Here ${_1F_1}$ denotes the standard hypergeometric function, see more detail in section 2.

\medskip
Our proof is based on elementary observation, with some interesting features. For example, to derive \eqref{staHUP0-DNbis} from \eqref{staHUP0}, the analysis in \cite{DLLZ} used the closeness of $\Sigma_0$ with respect to the semi-norm $\|\nabla v\|_2$; the attainability of $\inf_{u^* \in \Sigma_0}\|\nabla(u-u^*)\|_2$, and involved discussion according to the value of $\delta_1(u)$. Indeed, \eqref{staHUP0-DNbis} is a direct consequence of \eqref{staHUP0}, seeing the following abstract result in very general frameworks.

\begin{lem}\label{2.4}
Let $(\E, {\langle\cdot,\cdot\rangle})$ be an inner product space. Let $\Sigma \subset \E$ be a cone, that is $0 \in \Sigma$, and $\Sigma = \lambda\Sigma$ for all $\lambda \ne 0$. Given any $u \in \E$, we denote $S^u = \{w \in \Sigma, \|w\| = \|u\|\}$. There holds
\begin{align}
\label{d1}
{\rm dist}(u, S^u) \leq \sqrt{2}{\rm dist}(u, \Sigma).
\end{align}
Moreover, if $u \perp \Sigma$, then the equality holds true, i.e.~${\rm dist}(u, S^u) = \sqrt{2}{\rm dist}(u, \Sigma)$.
\end{lem}

Notice that the estimate \eqref{d1} does not require the completeness of $\E$, nor the closeness of $\Sigma$ in any topology, nor the attainability of the distance ${\rm dist}(u, \Sigma)$.

\smallskip
Lemma \ref{2.4} and some other preliminary results will be exhibited in section 2. The proof of Theorems \ref{1.1} and \ref{HH2-3} will be given in section 3. In section 4, we show the sharp stability estimate for the inequality \eqref{HUP0}, i.e.~Theorem \ref{sharpsta}.

\section{Preliminaries}
\reset
In this section, we state some useful Lemmas. The spherical harmonics decomposition is a powerful and important technic in analysis, in particular for the study of functional inequalities, some classical applications can be found in \cite{L, B, FLS, VZ}. Let $u \in C_c^\infty(\R^N)$, it is well known that $u$ can be expanded as follows
\begin{align}\label{HE}
u(x) = \sum_{k =0}^\infty r^kv_k(r)\phi_k(\sigma)
\end{align}
with $v_k\in C_c^2(\R_+)$, $\sigma = \frac{x}{|x|}$, and suitably normalized spherical eigenfunctions $\phi_k$ satisfying $$-\Delta_{S^{N-1}}\phi_k = k(N+k-2)\phi_k, \quad \forall\; k\in \N.$$ By Lemma 3.1 in \cite{CFL}, the following identities hold.
\begin{lem}\label{2.1}
Let $N\geq 2$, for any $u\in C_c^\infty(\R^N)$ and the sequence $v_k$ given by \eqref{HE}, we have
\begin{align*}
\int_{\R^N}|\nabla u|^2 dx =\sum_{k=0}^\infty \int_0^\infty r^{N+2k-1} |v'_k(r)|^2 dr,
\end{align*}
\begin{align*}
\int_{\R^N}\frac{|\nabla u|^2}{|x|} dx =\sum_{k=0}^\infty \Big(\int_0^\infty r^{N+2k-2} |v'_k(r)|^2 dr + k \int_0^\infty r^{N+2k-4} |v_k(r)|^2 dr\Big),
\end{align*}
\begin{align*}
\int_{\R^N}|x|^2|\nabla u|^2 dx =\sum_{k=0}^\infty\Big(\int_0^\infty r^{N+2k+1} |v'_k(r)|^2 dr - 2k \int_0^\infty r^{N+2k-1} |v_k(r)|^2 dr\Big)
\end{align*}
and
\begin{align*}
\int_{\R^N} |\Delta u|^2 dx =\sum_{k=0}^\infty \Big(\int_0^\infty r^{N+2k-1} |v''_k(r)|^2 dr + (N+2k-1) \int_0^\infty r^{N+2k-3} |v'_k(r)|^2 dr\Big).
\end{align*}
\end{lem}

By the above formulae, we observe that the spherical harmonic decomposition fit well with the understanding of HUP in additional forms as \eqref{HN+}, or the stability problem \eqref{staHUP+}.

\medskip
As already mentioned, another important tool for us is the Hardy type inequalities. Associated to the decomposition \eqref{HE}, we need especially the Hardy inequalities in $\R_+$ with general weights as follows, see for example \cite{GM}. 
A very general version without any symmetry assumption can be found in \cite{HY}.
\begin{lem}Let $V\in C^1(\R_+)$ be nonnegative, $f\in C^2(\R_+)$ be positive and $v\in C_c^1(\R_+)$ satisfy $V(0)v(0) = 0$. Then
\begin{align}
\label{Hardy-V}
\int_0^\infty V(r) |v'(r)|^2 dr \geq \int_0^\infty W(r) v^2(r) dr, \quad \mbox{with }\; W(r) = -\frac{(V f')'}{f}.
\end{align}
\end{lem}

The proof is the same as for classical Hardy inequality, by expanding
\begin{align}
\label{Hardy1}
\int_0^\infty V(r)\Big|v' - \frac{f'}{f}v\Big|^2(r) dr \geq 0,
\end{align}
and using integration by parts, so we omit the detail.

\smallskip
For example, we have a slight generalization of the classical Hardy inequalities.
\begin{align}\label{Hardy2}
\int_0^\infty r^\theta |v'(r)|^2 dr \geq \frac{(\theta-1)^2}{4}\int_0^\infty r^{\theta-2} v^2(r) dr, \quad \forall \;\theta > 1, \; v \in C_c^1(\R_+).
\end{align}
We can get \eqref{Hardy2} by applying formally \eqref{Hardy1} with $V = r^\theta$ and $f(r) = r^{-\frac{1-\theta}{2}}$. Notice that $v(0) = 0$ is not required. Here $f$ has a singularity at $r = 0$, but $\frac{f'}{f} = Cr^{-1}$, the corresponding weight $W \in L^1_{loc}(\R_+)$ and all computations with integration by parts for \eqref{Hardy1} work still.

We will also use \eqref{Hardy-V} with $V = r^\theta e^{-2r}$ and $f =e^{\kappa r}$, which yields
\begin{align*}
W:=-\frac{(V f')'}{f}=-\left[\theta\kappa r^{\theta-1} + \kappa(\kappa-2)r^\theta\right] e^{-2r}.
\end{align*}
\begin{lem}
\label{newHardy}
For any $\theta \ge 1$, $\kappa\in\R$ and $v\in C_c^1(\R_+)$, there holds
\begin{align}\label{HY}
\int_0^\infty r^\theta e^{-2r} |v'(r)|^2 dr \geq -\theta\kappa\int_0^\infty r^{\theta-1}e^{-2r} v^2(r) dr - \kappa(\kappa-2)\int_0^\infty r^\theta e^{-2r} v^2(r) dr.
\end{align}
\end{lem}

Next, we recall some elementary fact for hypergeometric functions, see \cite{V}. Denote
$$ H(t):= {_1F_1}\left(\beta;\alpha;t\right) = \sum_{j=0}^\infty \frac{\Gamma(j+\beta)}{\Gamma(j+\alpha)} \frac{t^j}{j!}, \quad \alpha, \beta \ge 0.$$
Then $H \in C^\infty(\R)$ and resolves
\begin{align}\label{HG}
tH''(t)+(\alpha-t)H'(t)-\beta H(t)=0 \quad \mbox{in } \R.
\end{align}
Our simple observation is
\begin{lem}
\label{sharpG}
Let $\Psi(t) = {_1F_1}\big(\beta;\alpha;-\frac{t^2}{2}\big)$ with $\alpha, \beta \ge 0$, then
\begin{align}
\label{GE}
t\Psi''(t)+t^2\Psi'(t)+(2\alpha-1)\Psi'(t)+ 2\beta t\Psi(t)=0 \quad \mbox{in } \R.
\end{align}
Moreover, for any $\ell \in\mathbb{N}$, there holds
\begin{align}\label{GF}
\Psi^{(\ell)}(t)=O\big(t^{-2\beta-\ell}\big),\quad \text{as } |t|\to\infty.
\end{align}
\end{lem}

\noindent
{\sl Proof.}
Clearly, $\Psi$ is pair and smooth in $\R$. Moreover
\begin{align*}
& \quad t\Psi''(t)+ t^2 \Psi'(t) + (2\alpha-1)\Psi'(t)+2\beta t\Psi(t)\\
& = t\Big[t^2H''(t)-H'(t)-t^2H'(t)-(2\alpha-1)H'(t)+2\beta H(t)\Big]\\
& = -2t\Big[tH''(t)+(\alpha-t)H'(t)-\beta H(t)\Big]
\end{align*}
So \eqref{GE} is satisfied seeing \eqref{HG}. Furthermore, when $t\to -\infty$, it is known \cite{V} that either ${_1F_1}\left(\beta;\alpha;t\right) = e^t$ if $\beta=\alpha$; and otherwise ${_1F_1}\left(\beta;\alpha;t\right) \sim C_{\beta,\alpha} |t|^{-\beta}$. Notice also that $H^{(\ell)}(t)={_1F_1}\left(\ell+\beta; \ell+\alpha;t\right)$ for any $\ell\in \mathbb{N}$. We get easily \eqref{GF} by induction. \qed

\subsection{Proof of Lemma \ref{2.4}}
Here we show the estimate \eqref{d1} which is valid in very general setting for stability problem. To simplify, for $u \in \E$ and $A \subset \E$, we denote the distance function in $\E$ as ${\rm d}(u, A) = \inf_{v\in A}\|u-v\|$.

\smallskip
If $u = 0$, \eqref{d1} is obviously true since $0 \in \Sigma\cap\Sigma^u$. Let $u \ne 0$, by scaling we can assume that $\|u\| = 1$, hence ${\rm d}(u, \Sigma) \in {[0, 1]}$.
\begin{itemize}
\item[Case 0.] If now $u \perp \Sigma$, clearly ${\rm d}(u, \Sigma) = \|u\|$ and ${\rm d}(u, S^u) = \sqrt{2}\|u\|$, because $\|u - p\|^2 = \|u\| + \|p\|^2$, for any $p \in \Sigma$.
\item[Case 1.] Let ${\rm d}(u, \Sigma) =0$. There is $\{p_i\} \subset \Sigma$ such that $\|p_i - u\| \to 0$, so $\|p_i\| \to 1$. Let $w_i = \frac{p_i}{\|p_i\|}$, then $w_i \in S^u$ and $\|w_i - u\| \to 0$, hence ${\rm d}(u, S^u) = 0$.
\item[Case 2.] Let ${\rm d}(u, \Sigma) = 1$. Using $\|u - \lambda w\|^2 \geq 1$ for any $\lambda$, there holds $\langle u, w\rangle =0$, hence $u \perp \Sigma$, so we are in Case 0.
\item[Case 3.] Let ${\rm d}(u, \Sigma) = \gamma \in {(0, 1)}$. There is $\{p_i\} \subset \Sigma$ such that $\|p_i - u\| \to \gamma$. Denote by $q_i$ the orthogonal projection of $u$ over ${\rm span}\{p_i\}$, we see that
\begin{align*}
u-q_i \perp q_i, \;\mbox{so } \|u-q_i\|^2 = 1 - \|q_i\|^2 \quad\mbox{and}\quad \|u-q_i\| \to \gamma.
\end{align*}
Hence we can assume $0 < \|q_i\| < 1$ and define $w_i = \frac{q_i}{\|q_i\|} \in S^u$. There holds
  \begin{align*}
\|u-w_i\|^2 = 2 - 2\frac{\langle u, q_i\rangle}{\|q_i\|} = 2 - 2\|q_i\| \leq 2 - 2\|q_i\|^2 = 2\|u-q_i\|^2.
\end{align*}
We get then ${\rm d}(u, S^u) \leq \sqrt{2}{\rm d}(u, \Sigma)$.
\end{itemize}
To conclude, \eqref{d1} is always valid. \qed

\section{Hydrogen uncertainty principle}
\reset
Here we consider the best constants $\HH_N$ defined by \eqref{HN}, and prove Theorems \ref{1.1}-\ref{HH2-3}. Let $u \in C_c^\infty(\R^N)$ and
\begin{align}
\label{IN}
J_N(u) = \int_{\R^N} (\Delta u)^2 dx + \int_{\R^N} |\nabla u|^2 dx - (N+1)\int_{\R^N} \frac{|\nabla u|^2}{|x|} dx
\end{align}
Using the notation of spherical harmonic expansion \eqref{HE} and Lemma \ref{2.1}, we have
\begin{align}
\label{ineq1}
J_N(u) = \sum_{k=0}^\infty J_{N,k}(v_k)
\end{align}
where
\begin{align}
 \label{JNK}
\begin{split}
J_{N, k}(v) & = \int_0^\infty r^{N+2k-1}v''(r)^2 dr + (N+2k-1)\int_0^\infty r^{N+2k-3}v'(r)^2 dr\\
& \quad + \int_0^\infty r^{N+2k-1}v'(r)^2 dr - (N+1) \int_0^\infty r^{N+2k-2}v'(r)^2 dr\\
& \quad - (N+1)k\int_0^\infty r^{N+2k-4}v(r)^2 dr.
\end{split}
\end{align}
We will estimate each term $J_{N, k}(v)$ separately. 
\begin{lem}
\label{easyJ}
Let $v \in C_c^2(\R_+)$. Then $J_{N, 0}(v) \geq 0$ for all $N \ge 1$. Moreover, $J_{N, k}(v) \geq 0$ for all $k \geq 2$ and $N \geq 2$.
\end{lem}

\noindent
{\sl Proof.} Direct calculations yield
\begin{align}
 \label{JNKv}
\begin{split}
J_{N, k}(v) & = \int_0^\infty r^{N+2k-1}(v''+v')^2 dr + (N+2k-1)\int_0^\infty r^{N+2k-3}v'(r)^2 dr \\
& \quad + 2(k-1) \int_0^\infty r^{N+2k-2}v'(r)^2 dr - (N+1)k\int_0^\infty r^{N+2k-4}v(r)^2 dr.
\end{split}
\end{align}

First, we consider $J_{N, 0}(v)$. Let $v = e^{-r}\xi(r)$, there holds $v'' + v' = (v'+v)' = e^{-r}(\xi''-\xi')$. Denote $\rho = \xi'-\xi$, then
\begin{align}
 \label{ineq2.0}
\begin{split}
J_{N, 0}(v) & = \int_0^\infty r^{N-1}(v''+v')^2 dr + (N-1)\int_0^\infty r^{N-3}v'(r)^2 dr -2 \int_0^\infty r^{N-2}v'(r)^2 dr\\
& = \int_0^\infty r^{N-1}e^{-2r}\rho'(r)^2 dr + (N-1)\int_0^\infty r^{N-3}e^{-2r}\rho(r)^2 dr\\
& \quad - 2 \int_0^\infty r^{N-2}e^{-2r}\rho(r)^2 dr\\
& = \int_0^\infty r^{N-3}e^{-2r} \big(r\rho' - \rho\big)^2 dr\geq 0.
\end{split}
\end{align}

Let $k \geq 2$ and $N \geq 3$. By \eqref{Hardy2} with
$\theta = N+2k-2$, we get
\begin{align}
 \label{ineq2.k.2}
\begin{split}
2(k-1) \int_0^\infty r^{N+2k-2}v'(r)^2 dr & \geq \frac{(k-1)(N+2k-3)^2}{2}\int_0^\infty r^{N+2k-4}v(r)^2 dr\\
& \geq (N+1)k\int_0^\infty r^{N+2k-4}v(r)^2 dr.
\end{split}
\end{align}
The last line holds since $2(k-1)(N+2k-3)^2\geq k(N+1)^2\geq 4k(N+1)$. Inserting \eqref{ineq2.k.2} into \eqref{JNKv}, we have $J_{N, k}(v) \geq 0$ for $k \geq 2$, $N\geq 3$.

\medskip
It remains to handle $J_{2,2}(v)$. Let again $v=e^{-r}\xi(r)$, then $v''+v'= (e^{-r}\xi')'$. Moreover, applying \eqref{HY} with $(\theta, \kappa) = \big(3, \frac{1}{3}\big)$, and $(\theta, \kappa) = (4,1)$ respectively: For any $\xi \in C_c^1(\R_+)$, there hold
\begin{align*}
\int_0^\infty r^{3}e^{-2r}\xi'^2 dr &\geq \frac{5}{9}\int_0^\infty r^{3} e^{-2r}\xi^2 dr -\int_0^\infty r^{2}e^{-2r}\xi^2 dr,\\
\int_0^\infty r^{4}e^{-2r}\xi'^2 dr &\geq \int_0^\infty r^{4} e^{-2r}\xi^2 dr -4\int_0^\infty r^{3}e^{-2r}v\xi^2 dr.
\end{align*}
Combining \eqref{JNKv}, \eqref{Hardy2}, integration by parts and the above two inequalities, we see that
\begin{align*}
J_{2,2}(v) & \geq 4\int_0^\infty r^{3} e^{-2r}\xi'^2 dr +5\int_0^\infty r^{3}e^{-2r}(\xi'-\xi)^2 dr\\
& \quad +2\int_0^\infty r^{4}e^{-2r}(\xi'-\xi)^2 dr -6\int_0^\infty r^{2} e^{-2r}\xi^2 dr\\
& = 9\int_0^\infty r^{3}e^{-2r}\xi'^2 dr +2\int_0^\infty r^{4}e^{-2r} \xi'^2 dr -2\int_0^\infty r^{4} e^{-2r} \xi^2 dr \\
&\quad +3\int_0^\infty r^{3} e^{-2r}\xi^2 dr +9\int_0^\infty r^{2}e^{-2r} \xi^2 dr\\
& \geq 0.
\end{align*}
The proof is completed. \qed

\subsection{Sharp inequality for $N \geq 4$}
We will show Theorem \ref{1.1} with \eqref{HN+} and
\begin{prop}
Let $N\geq 4$. Then $J_{N,k}(v)\geq 0$ for any $k\geq 0$ and $v \in C_c^2(\R_+)$. Hence $J_N(u) \geq 0$ for any $u \in W^{2,2}(\R^N)$. Moreover, $J_N(u) =0$ if and only if $u(x) = \alpha(1+r)e^{-r}$, $\alpha \in \R$.
\end{prop}

\smallskip
\noindent
{\sl Proof.} Thanks to Lemma \ref{easyJ}, we focus only on $J_{N,1}(v)$. Let $v \in C_c^2(\R^+)$, set $v = e^{-r}\xi(r)$, so $v'' + v' = \big(e^{-r}\xi')'$. Using again the Hardy inequality \eqref{Hardy2}, as $N\geq 4$
\begin{align}
 \label{ineq2}
\begin{split}
J_{N, 1}(v) & \geq \frac{N^2}{4}\int_0^\infty r^{N-1}e^{-2r}\xi'(r)^2 dr + (N+1)\int_0^\infty r^{N-1}e^{-2r}(\xi' - \xi)^2 dr\\
&\quad - (N+1)\int_0^\infty r^{N-2}e^{-2r}\xi(r)^2 dr\\
& = \frac{(N+2)^2}{4}\int_0^\infty r^{N-1}e^{-2r}\xi'(r)^2 dr + (N+1)(N-2)\int_0^\infty r^{N-2}e^{-2r}\xi(r)^2 dr\\
&\quad  - (N+1)\int_0^\infty r^{N-1}e^{-2r}\xi(r)^2 dr\\
& =: R_{N,1}(\xi).
\end{split}
\end{align}

Apply the Hardy inequality \eqref{HY} with $\theta = {N-1}$ and $\kappa = \frac{2}{N+2}$: For all $\xi \in C_c^1(\R_+)$, there holds
\begin{align}
 \label{Hardy3}
 \begin{split}
\int_0^\infty r^{N-1}e^{-2r}\xi'(r)^2 dr & \geq \frac{4(N+1)}{(N+2)^2}\int_0^\infty r^{N-1}e^{-2r}\xi(r)^2 dr\\
& \quad - \frac{2(N-1)}{N+2}\int_0^\infty r^{N-2}e^{-2r}\xi(r)^2 dr.
\end{split}
\end{align}
Inserting \eqref{Hardy3} into the expression of $R_{N,1}$, we obtain
\begin{align*}
R_{N,1}(\xi) &\geq \Big[(N+1)(N-2) - \frac{(N-1)(N+2)}{2}\Big]\int_0^\infty r^{N-2}e^{-2r}\xi(r)^2 dr\\
& = \frac{N^2 -3N-2}{2}\int_0^\infty r^{N-2}e^{-2r}\xi(r)^2 dr,
\end{align*}
hence $J_{N.1}(v) \geq R_{N,1}(\xi) \geq 0$ if $N \geq 4$. This means $J_N(u) \geq 0$ for any $u \in C_c^2(\R^N)$ and $N \geq 4$. Therefore \eqref{CFL} holds true in $W^{2,2}(\R^N)$ by density argument.

\medskip

Moreover, by the above analysis, noticing the non-attainability of Hardy's inequalities \eqref{Hardy1} for any nontrivial function, we see that to reach the best constant $(N+1)$ in \eqref{CFL}, the only choice is $v_k \equiv 0$ for $k \geq 1$; and the equality holds in \eqref{ineq2.0}. This means that $r\rho' \equiv \rho$, that is $\xi'-\xi = \rho = Cr$. Finally $u(x) = v_0(r) = e^{-r}\xi \in W^{2,2}(\R^N)$ yields that $v_0(r) = -C(1+r)e^{-r}$ with $C \in \R$. \qed

\begin{rem}
The cone of extremal functions for $J_N$ loses one degree of freedom comparing to \eqref{CFL}, since $J_N$ is no longer scaling invariant, conversely to \eqref{CFL}.
\end{rem}

\subsection{Estimates for $N=2$ and $3$}
Recall that $\HH_N^+ \leq (N+1)$ since $J_N(u_0) = 0$ with $u_0(x) = (1+r)e^{-r}$. Combining with Lemma \ref{easyJ}, we know that $\HH_N^+ = \min(N+1, \CC_N)$ with
\begin{align}
\label{CCN}
\CC_N := \inf_{v \in C_c^2(\R_+), v \ne 0}\frac{E_1(v)}{F_1(v)},
\end{align}
where
\begin{align*}
E_1(v) := \int_0^\infty r^{N+1}v''(r)^2 dr + (N+1)\int_0^\infty r^{N-1}v'(r)^2 dr + \int_0^\infty r^{N+1}v'(r)^2 dr
\end{align*}
and
\begin{align*}
F_1(v) := \int_0^\infty r^Nv'(r)^2 dr + \int_0^\infty r^{N-2}v(r)^2 dr.
\end{align*}
Chen-Tang's counter example $u_1 = re^{_r}\phi_1(\sigma)$ means that $\HH_N^+ = \CC_N$ for $N=2$ or $3$. So we consider $G_\a(v) = E_1(v)- \a F_1(v)$ with $\alpha > 0$, and we are trying to estimate the constant $\a$ such that $G_\a(v) \geq 0$ in $C_c^2(\R_+)$.

\medskip
$\bullet$ Consider firstly $N= 2$. Then, for any $\epsilon \leq 1$,
\begin{align}
\label{G2a}
\begin{split}
G_{\alpha}(v) & = \int_0^\infty  r^3 v''^2 dr + 3\int_0^\infty r v'^2dr + \int_0^\infty r^3 v'^2 dr -\alpha \int_0^\infty r^2 v'^2 dr -\alpha \int_0^\infty v^2 dr\\
& \geq \big(2\sqrt{1-\epsilon}+1-\alpha\big)\int_0^\infty r^2 v'^2 dr +3 \int_0^\infty r v'^2 dr +\epsilon \int_0^\infty r v^2 -\alpha \int_0^\infty v^2 dr.
\end{split}
\end{align}
Here we used the following inequality given by Hamamoto \cite{H} (see also \cite{HOT}) with $\mu=2$: Given $\mu, \epsilon\in\R$ and $\epsilon\leq \frac{\mu^2}{4}$, there holds, for any $f\in C_c^2(\R_+)$,
\begin{align*}
& \quad \int_0^\infty x^{\mu+1} f''(x)^2 dx + \int_0^\infty x^{\mu+1} f'(x)^2 dx\\
& \geq \epsilon \int_0^\infty x^{\mu-1} f(x)^2 dx +\big(\sqrt{\mu^2-4\epsilon}+1\big)\int_0^\infty x^\mu f'(x)^2 dx.
\end{align*}
Rewrite again $v = e^{-r}\xi(r)$, using integration by parts, there holds
\begin{align}
\label{G2b}
\begin{split}
\int_0^\infty r v'^2 dr & = \int_0^\infty r e^{-2r} \xi'^2 dr - \int_0^\infty r e^{-2r} \xi^2 dr +\int_0^\infty e^{-2r} \xi^2 dr\\
& \geq \frac{\alpha}{3}\int_0^\infty e^{-2r} \xi^2 dr - \frac{\alpha^2}{9} \int_0^\infty r e^{-2r} \xi^2 dr.
\end{split}
\end{align}
For the last line, we applied the Hardy inequality \eqref{HY} with $(\theta, \kappa) = \big(1, {1-\frac{\alpha}{3}}\big)$.

Combining \eqref{G2a} and \eqref{G2b}, we arrive at
\begin{align*}
G_\alpha(v)  \geq \big(2\sqrt{1-\epsilon}+1-\alpha\big)\int_0^\infty r^2 v'^2 dr +\Big(\epsilon-\frac{\alpha^2}{3}\Big) \int_0^\infty r v^2.
\end{align*}
By Chen-Tang's estimate, we need only to consider $\a \leq \sqrt 3$. Fix $\epsilon =\frac{\alpha^2}{3}\leq 1$, then $G_\alpha(v) \geq 0$ if $$2\sqrt{1-\epsilon}+1-\alpha = 2\sqrt{1-\frac{\alpha^2}{3}}+1-\alpha\geq 0.$$ The maximum value satisfying the above inequality is $\a = \frac{3+6\sqrt 2}{7}$, i.e.~$\HH_2^+ \geq \frac{3+6\sqrt 2}{7} \simeq 1.6407$.

\medskip
$\bullet$ Consider now $N=3$. Let $G_\a$, $E_1$, $F_1$ be as above and $\alpha > 0$. Let $v = e^{-r}\xi(r)$, so $v''+v'=(e^{-r}\xi')'$. By the Hardy inequality \eqref{Hardy2} with $\theta = 4$ and integration by parts, we get
\begin{align*}
G_4(v) & = \int_0^\infty r^4(v''+v')^2 dr + 4\int_0^\infty r^2v'(r)^2 dr - 4\int_0^\infty rv(r)^2 dr\\
&\geq \frac{9}{4}\int_0^\infty r^2 e^{-2r}\xi'(r)^2 dr + 4\int_0^\infty r^2 e^{-2r}(\xi'-\xi)^2(r) - 4\int_0^\infty re^{-2r}\xi^2(r) dr\\
& =\frac{25}{4}\int_0^\infty r^2e^{-2r}\xi'(r)^2 dr + 4\int_0^\infty re^{-2r}\xi(r)^2 dr - 4\int_0^\infty r^2e^{-2r}\xi(r)^2 dr.
\end{align*}
Therefore, for $\a \leq 4$ and $\xi \in C_c^1(0, \infty)$,
\begin{align*}
G_\alpha(v) & = G_4(v) + (4-\a)F_1(v)\\
& \geq \frac{25}{4}\int_0^\infty r^2e^{-2r}\xi'(r)^2 dr + (8-\alpha)\int_0^\infty re^{-2r}\xi(r)^2 dr - 4\int_0^\infty r^2e^{-2r}\xi(r)^2 dr\\
&\quad + (4-\a)\int_0^\infty r^3\big(e^{-r}\xi)'(r)^2 dr\\
& \geq \frac{25}{4}\int_0^\infty r^2e^{-2r}\xi'(r)^2 dr - 4\int_0^\infty r^2e^{-2r}\xi(r)^2 dr + (12-2\alpha)\int_0^\infty re^{-2r}\xi(r)^2 dr.
\end{align*}
In the last line, we have lower bounded the integral of $r^3\big(e^{-r}\xi)'(r)^2$ by \eqref{Hardy2} with $\theta = 3$.

Moreover, applying \eqref{HY} with $(\theta, \kappa) = \big(2, \frac{2}{5}\big)$, there holds
\begin{align}
 \label{ineq3}
\int_0^\infty r^2e^{-2r}\xi'(r)^2 dr \geq \frac{16}{25}\int_0^\infty r^2e^{-2r}\xi(r)^2 dr - \frac{4}{5}\int_0^\infty re^{-2r}\xi(r)^2 dr.
\end{align}
We conclude then
$$G_\alpha(v) \geq (7-2\a)\int_0^\infty re^{-2r}\xi(r)^2 dr, \quad \forall\; v\in C_c^2(\R_+),$$
hence $\HH_3^+ \geq \frac{7}{2}$.

\medskip
Using \eqref{HN+}, we get readily \eqref{newHN2-3}. Even we cannot fix the values of $\HH_2 = \CC_2$ and $\HH_3 = \CC_3$, an interesting observation is that they are reached by some extremal functions. Indeed, we will show that $\CC_2$ and $\CC_3$ admit minimizers. We only explain for $N=2$, the three dimensional case is just similar. Define $\|v\|_* = \|\Delta v\|_2 + \|\nabla v\|_2$, and
$$\M_1 = \big\{u(x) = rv(r)\phi_1(\sigma)\big\}\cap \overline{C_c^\infty(\R^2)}^{\|\cdot\|_*}$$
Here $\phi_1(\sigma)$ means all normalized function with linear combination of $\cos\sigma$ and $\sin\sigma$, so the set of $\phi_1(\sigma)$ is homeomorphic to $S^1$, hence compact. By Lemma \ref{2.1}, $u(x) = rv(r)\phi_1(\sigma) \in \M_1$ if and only if $E_1(v) < \infty$, therefore, to show that $\CC_2$ is attained in $\M_1$, we need only to prouve the compactness of the mapping $v \mapsto F_1(v)$ over $\M_1$.

Clearly, by Lemma \ref{2.1}, $\|(r + r^3)v'^2\|_{L^1(\R_+)} \leq E_1(v)$, and $\sqrt{E_1(v)}$ defines a norm equivalent to $\|v\|_*$ in $\M_1$. So we get uniform control for the integral of $r^2v'(r)^2$ near $0$ and $+\infty$ by $E_1(v)$;. On the other hand, with classical Sobolev embedding, the mapping $v \mapsto \|rv'\|_{L^2(I)}$ is compact over any interval $I \subset\subset (0,\infty)$ with respect to $E_1(v)$, so is the mapping $v \mapsto \|rv'\|^2_{L^2(\R_+)}$.

By $G_1(v) \geq 0$ and the Hardy inequality \eqref{Hardy2}, there holds $\|v\|_{L^2(\R_+)}^2 + \|rv^2\|_{L^1(\R_+)} \leq C$ if $E_1(v)$ is bounded. We have then uniform control to the integral of $v^2$ near $+\infty$. Moveover, take $\beta \in (1,2)$, as $\|r^\beta v'^2\|_{L^1(\R_+)}$ is uniformly bounded,
$$|v(r)|^2 \leq C\int_r^\infty s^{-\beta}ds \leq Cr^{1-\beta}.$$
This yields the equi-integrability of $v^2$ near $r=0$. As before, the integral of $v^2$ over compact intervals in $(0, \infty)$ can be handled by classical Sobolev embeddings. Finally the mapping $F_1$ is compact with respect to $E_1(u)$ over $\M_1$, thus $\CC_2$ is reached.
\qed

\section{Sharp stability for a Heisenberg Uncertainty Principles}
\reset
Here we consider the stability problem associated to the sharp inequality \eqref{HUP0}. Recall that the equality holds in \eqref{HUP0} if and only if $u \in \Sigma_0:= \big\{\alpha e^{-\beta r^2}, \; \alpha\in \R, \; \beta>0 \big\}.$ Denote
\begin{align*}
H_0(u) = \|\D u\|_2^2, \quad U_0(u) = \|r\nabla u\|_2^2, \quad P_0(u) = \|\nabla u\|_2^2.
\end{align*}
Let $\H_0$ be the completion of $C_c^\infty(\R^N)$ under the norm $\|\D u\|_2 + \|r\nabla u\|_2$. Our aim is to find the constant
\begin{align}
\label{staHUP0.1}
\SS_0 = \inf_{u \in \H_0, d_0(u, \Sigma_0) > 0} \frac{\sqrt{H_0(u)U_0(u)} - \frac{N+2}{2}P_0(u)}{d_0(u, \Sigma_0)^2}
\end{align}
where 
$$d_0(u, \Sigma_0) = \inf_{v^* \in \Sigma_0} \|\nabla (u-v^*)\|_2.$$
Thanks to Lemma \ref{newlem1}, we need only to determine
\begin{align}
\label{staHUP0.1+}
\SS_{0, +} = \inf_{u \in \H_0, d_0(u, \Sigma_0)> 0} \frac{\delta_2(u)}{d_0(u, \Sigma_0)^2} = 2\SS_0
\end{align}
with
\begin{align*}
\delta_2(u):=\int_{\R^N} |\Delta u|^2 dx + \int_{\R^N} |x|^2|\nabla u|^2 dx - (N+2)\int_{\R^N} |\nabla u|^2 dx.
\end{align*}

Now we recall briefly some analysis in \cite{DLLZ}. Let $u \in C_c^\infty(\R_+)$, using the harmonic expansion \eqref{HE} of $u$ and Lemma \ref{2.1}, we have
\begin{align*}
\delta_2(u) = \sum_{k=0}^\infty I_{N,k}(v_k),
\end{align*}
where for all $N\geq 2$, $k \in \N$, $v\in C_c^2(\R_+)$,
\begin{align*}
I_{N,k}(v) & := \int_0^\infty r^{N+2k-1}|v''(r)|^2 dr + (N+2k-1)\int_0^\infty r^{N+2k-3}|v'(r)|^2 dr\\
& \quad +\int_0^\infty r^{N+2k+1}|v'(r)|^2 dr -2k\int_0^\infty r^{N+2k-1}|v(r)|^2 dr\\
& \quad - (N+2)\int_0^\infty r^{N+2k-1}|v'(r)|^2 dr.
\end{align*}
Denote $u_k = r^kv_k(r)\phi_k(\sigma)$ the general terms in the expansion \eqref{HE} and
\begin{align}
\label{inner}
\langle u, w\rangle = \int_{\R^N} \nabla u\cdot\nabla w dx.
\end{align}
Recall that $u_k \perp u_\ell$ for $k\ne \ell$, and $u_k \perp \Sigma_0$ for $k \geq 1$. Consequently,
\begin{align*}
d_0(u, \Sigma_0)^2 & = \sum_{k=0}^\infty d_0(u_k, \Sigma_0)^2 = d_0(u_0, \Sigma_0)^2 + \sum_{k=1}^\infty \|\nabla u_k\|_2^2.
\end{align*}
Moreover, by Lemma \ref{2.1}, for any $k \in \N$,
\begin{align*}
\|\nabla u_k\|_2^2 = Q_k(v_k), \quad \mbox{where } Q_k(v) := \int_0^\infty r^{N+2k-1} |v'(r)|^2 dr.
\end{align*}
Define
\begin{align*}
\CC_{N,0} = \inf_{u \in C_{c, {\rm rad}}^\infty(\R^N)}\frac{\delta_2(u)}{d_0(u, \Sigma_0)^2}, \quad \mbox{and} \quad\forall\; k \geq 1,\; \CC_{N,k} = \inf_{v \in C_c^2(\R_+)}\frac{I_{N,k}(v)}{Q_k(v)}.
\end{align*}
Do-Lam-Lu-Zhang (see \cite[proof of Theorems 1.2-1.3]{DLLZ}) showed that $\CC_{N,0}\geq 2$ and for $k \geq 1$,
\begin{align*}
\sqrt{(N+2k)^2 - 8k} - N \leq \CC_{N,k} \leq 2k, \quad \forall\; k \geq 1.
\end{align*}
Therefore
$$\SS_{0, +} = \inf_{k \in \N} \CC_{N, k} \geq \min_{k \geq 1}\sqrt{(N+2k)^2 - 8k} - N = \sqrt{N^2 +4N-4} - N.$$

Our key observation is the following
\begin{lem}
For any $k \geq 1$, $\CC_{N, k} = \sqrt{(N+2k)^2 - 8k} - N$. Hence
\begin{align}
\label{SS0+}
\SS_{0, +} = \CC_{N,1} = \sqrt{N^2 +4N-4} - N.
\end{align}
\end{lem}

\noindent
{\sl Proof.} By the computation in \cite[page 18]{DLLZ}, for any $k \in \N$, $\gamma \in \R$, $v \in C_c^2(\R_+)$, we have the following equality
\begin{align*}
& I_{N,k}(v) - \gamma Q_k(v)\\
=&\int_0^\infty \Big|rv''(r)+ r^2 v'(r)+ A_k v'(r) + B_{k,\gamma} rv(r)\Big|^2 r^{N+2k-3} dr + T_k(\gamma)R_k(v)
\end{align*}
with
\begin{align*}
A_k = N + 2k-1, \quad B_{k,\gamma} = N+ k +\frac{\gamma}{2},\quad R_k(v):=\int_0^\infty r^{N+2k-1} v^2 dr
\end{align*}
and
\begin{align*}
T_k(\gamma) = (N+2k)B_{k,\gamma} - B_{k,\gamma}^2 -2k = (N-2+k)k - \frac{N}{2}\gamma - \frac{1}{4}\gamma^2.
\end{align*}
Clearly, $T_k(\gamma_k) = 0$ with $\gamma_k := \sqrt{(N+2k)^2 - 8k} - N > 0$; and $T_k(\gamma) < 0$ for $\gamma > \gamma_k$.

It means that $\CC_{N, k} \geq \gamma_k$. Moreover, we will have $I_{N,k}(v) - \gamma_k Q_k(v) \equiv 0$ if $v$ resolves the differential equation
$$rv'' + r^2 v' + A_k v' + B_{k,\gamma_k} rv = 0 \quad \mbox{in } \R_+.$$
By Lemma \ref{sharpG},
$$\widetilde v_k(r) := {_1F_1}\Big(\frac{B_{k, \gamma_k}}{2};\frac{A_k+1}{2};-\frac{r^2}{2}\Big)$$
is a such solution. Moreover, $\widetilde v_k \in C^\infty(\R_+)$ and by \eqref{GF}, $\widetilde v_k^{(\ell)}(r) = O(r^{-2N-2k-\gamma_k-\ell})$ as $r \to \infty$, for any $\ell \in \N$. So that $I_{N,k}(\widetilde v_k) + Q_k(\widetilde v_k) < \infty$, hence $\widetilde u_k(x) = r^k\widetilde v_k(r)\phi_k(\sigma)$ belongs to $\H_0$. For $\gamma > \gamma_k$ and $k \geq 1$, there holds
\begin{align*}
\delta_2(\widetilde u_k) - \gamma d_0(\widetilde u_k, \Sigma_0)^2 & = \delta_2(\widetilde u_k) - \gamma\|\nabla \widetilde u_k\|_2^2\\
& = I_{N,k}(\widetilde v_k) - \gamma Q_k(\widetilde v_k)\\
& = T_k(\gamma) R_k(\widetilde v_k) < 0.
\end{align*}
We deduce then $\CC_{N, k} \leq \gamma_k$ by standard approximation argument, so $\CC_{N, k} = \gamma_k$ for all $k \geq 1$, consequently \eqref{SS0+} holds true. \qed

\medskip
\noindent
{\sl Proof of Theorem \ref{sharpsta} completed.} By Lemma \ref{newlem1}, $\SS_{0, +} = 2\SS_0 = \CC_{N,1}$, so \eqref{staHUP0} is sharp.
Furthermore, as $\CC_{N, k} > \CC_{N,1}$ for $k\ne 1$, to reach the sharp equality $\delta_2(u) - \CC_{N,1} d_0(u, \Sigma_0)^2 = 0$, we must have $u = Cr\widetilde v_1(r)\phi_1(\sigma)$ with
$$\widetilde v_1(r)={_1F_1}\Big(\frac{N+1}{2} + \frac{\CC_{N,1}}{4};\frac{N+2}{2};-\frac{r^2}{2}\Big).$$ By scaling argument, we get all extremal functions to \eqref{staHUP0}.

\medskip
Now we consider \eqref{staHUP0-DNbis}. For any $u \in \H_0$, set
\begin{align*}
\SS^u_{0, +} = \inf_{u \in \H_0, d_0(u, S^u_0)> 0} \frac{\delta_2(u)}{d_0(u, S^u_0)^2} \quad\mbox{where } S^u_0 = \big\{w\in \Sigma_0, \|\nabla w\|_2 = \|\nabla u\|_2\big\}.
\end{align*}
 We can apply Lemma \ref{2.4} with $\E_0 = \H_0$ and the inner product given in \eqref{inner}. Therefore \eqref{d1} means that $\SS^u_{0, +} \geq \frac{\CC_{N,1}}{2}$. Take $\widetilde u_1$ as above, then $\widetilde u_1\perp \Sigma_0$ with respect to \eqref{inner}. By Lemma \ref{2.4}, there holds
$$\delta_2(\widetilde u_1) = \CC_{N,1} d_0( u, \Sigma_0)^2 = \frac{\CC_{N,1}}{2} d_0(u, S^u_0)^2.$$
Hence we get $\SS^u_{0, +} = \frac{\CC_{N,1}}{2}$. Finally, the sharpness of \eqref{staHUP0-DNbis} is ensured by the scaling argument as for \eqref{HN+} or Lemma \ref{newlem1}, so we are done.
\qed

\begin{rem}
 The above analysis is very practical in many other $L^2$ settings, for example to derive quickly \cite[Theorem 3.3]{CFLL}. We can also claim the stability of the sharp inequality in \eqref{CFL} or Theorem \ref{1.1} by similar approach, in other words, for $N\geq 4$, there is a constant $\SS_N > 0$, such that for all $u$ belonging to the closure of $C_c^\infty(\R^N)$ under the norm $\|\D v\|_2+ \|\nabla v\|_2$, there holds
\begin{align}
\label{CFLsta}
\|\Delta u\|_2 \|\nabla u\|_2 - \frac{N+1}{2}\int_{\R^N}\frac{|\nabla u|^2}{|x|} dx \geq \SS_N\inf_{u^* \in \Sigma_1}\int_{\R^N} \frac{|\nabla(u - u^*)|^2}{|x|} dx, 
\end{align}
here $\Sigma_1 = {\rm span}\{(1+r)e^{-r}\}$. However, we were not able to fix the best constant for \eqref{CFLsta}.
\end{rem}

\bigskip
\noindent
{\bf Acknowledgements.} The authors would like to thank Professor Guozhen Lu for showing us some precisions on \cite{DLLZ}. The authors are partially supported by NSFC (No.~12271164), Science and Technology Commission of Shanghai Municipality (No.~22DZ2229014). 

\smallskip
\noindent \textbf{Data availability.} Data sharing is not applicable to this work as no new data
were created or analyzed in the study.

\end{document}